\documentclass[12pt]{article}
\usepackage{amsmath,latexsym,amssymb,graphicx,epsf,amsfonts,amsthm}
\usepackage{epsf}
\usepackage{color}
\numberwithin{equation}{section} \numberwithin{figure}{section}

\begin{document}
\hoffset = -1truecm \voffset = -1truecm
\title{Computational Modeling of Spectral Data Fitting with Nonlinear Distortions}
\author{Yuanchang Sun\thanks{Department of Mathematics and Statistics,
Florida International University, Miami, FL 33199, USA.},  Wensong Wu$^{*}$, Jack  Xin\thanks{Department of Mathematics, University of California
at Irvine, Irvine, CA 92697, USA.} }
\date{}
\maketitle

\begin{abstract}
Substances such as chemical compounds are invisible to human eyes, they are usually captured by sensing equipments with their spectral fingerprints.  Though spectra of pure chemicals can be identified by visual inspection, the spectra of their mixtures take a variety of complicated forms.  Given the knowledge of spectral references of the constituent chemicals, the task of data fitting is to retrieve their weights, and this usually can be obtained by solving a least squares problem.  Complications occur if the basis functions (reference spectra) may not be used directly to best fit the data.  In fact, random distortions (spectral variability) such as shifting,  compression, and expansion have been observed in some source spectra when the underlying substances are mixed.  In this paper, we formulate mathematical model for such nonlinear effects and build them into data fitting algorithms.  If minimal knowledge of the distortions is available, a deterministic approach termed {\it augmented least squares} is developed and it fits the spectral references along with their derivatives to the mixtures.  If the distribution of the distortions is known a prior, we consider to solve the problem with maximum likelihood estimators which incorporate the shifts into the variance matrix.  The proposed methods are substantiated with numerical examples including data from Raman spectroscopy (RS), nuclear magnetic resonance (NMR), and differential optical absorption spectroscopy (DOAS) and show satisfactory results.

\end{abstract}
\section{Introduction}
Given a (or multiple) noisy spectral observation(s) of a sample containing several sources (e.g., chemicals), one attempts to retrieve the weights of source signals in the mixtures.  This can be solved by the least squares type of methods, which are the most widely used data analysis tools with numerous applications.  Suppose the observations are linear mixtures of the source signals
\begin{equation}
\label{LMM}
{\boldsymbol X} = \boldsymbol{A}\,\boldsymbol{S} + \boldsymbol{N}\;, \mathrm{with}\; \boldsymbol{X}_{ij} = x_{ij}, \boldsymbol {A}_{ij}= a_{ij}, \boldsymbol{S}_{ij} = s_{ij}, \boldsymbol{N}_{ij} = n_{ij}\;,
\end{equation}
where $\boldsymbol X \in \mathbb{R}^{m\times p} $ is the $m\times p$ observation matrix (mixture matrix), each row $\boldsymbol{x}_i$ stands for a single observation of length $p$; $\boldsymbol{S} $ is the basis (or source) matrix with each $\boldsymbol{s}_i$ being a pure source signal; $\boldsymbol{A} $ is the mixing matrix and has size $n\times m$, and $\boldsymbol{N} $ is the noise matrix of the same size as $\boldsymbol{X}$. Given $\boldsymbol{S}, \boldsymbol{X}$, the least squares solution of the mixing matrix is given by $\hat{\boldsymbol{A}} = \boldsymbol{X}\boldsymbol{S}^{\mathrm{T}} (\boldsymbol{S} \boldsymbol{S}^{\mathrm{T}})^{-1}$ which is the minimizer of the following problem
\begin{equation}
\label{LS}
\min_{\boldsymbol{A}}\|\boldsymbol{X}-\boldsymbol{A}\,\boldsymbol{S}\|^2_2\;.
\end{equation}
If the noise is normally distributed and variances of the observations are equal, the above ordinary least squares (OLS) method provides minimum-variance mean-unbiased estimation, and is actually the maximum likelihood estimator.  When the variances of the observations are unequal or certain correlations exist among the observations, the OLS can be inefficient and lead to poor estimates on the mixing matrix.  In some cases, one can apply a generalized least squares (GLS) model,
\begin{equation}
\label{GLS}
\boldsymbol{A}= \arg \min (\boldsymbol{X} - \boldsymbol{AS})\boldsymbol{Q}^{-1}(\boldsymbol{X}-\boldsymbol{SA})^{\boldsymbol{T}}
\end{equation}
with $\boldsymbol{Q}$ being the covariance matrix of the noise.  The estimate of $\boldsymbol{A}$ has an explicit formula
$\boldsymbol{A} = \boldsymbol{X} \boldsymbol{Q}^{-1}\boldsymbol{S}^{\mathrm{T}} (\boldsymbol{S}\boldsymbol{Q}^{-1}\boldsymbol{S}^{\mathrm{T}})^{-1}$
In practice, one needs to estimate $\boldsymbol{Q}$ for the implementation of GLS since the covariance matrix is in general unknown.  Hence an implementable version of GLS usually contains two steps; We first solve an OLS problem (or other variant least squares) to obtain an estimate of $\boldsymbol{Q}$ by the fitting residuals; Then we implement the GLS by (\ref{GLS}).  Least squares type methods have been very successful with numerous applications for data fitting.  However, complications occur if the basis functions (reference spectra $\boldsymbol S$) may not be used directly to best fit the data.  For example, the reference spectra may have either shifting or compression/expansion when the source objects are present in a mixture, these effects are also called spectral variability. Then the reference spectra need to be properly registered before fitting.

In this paper, we are concerned with a regime where the source signals possess such spectral variability in addition to the noise and measurement errors.  This has been a challenging problem in signal processing, for example, in hyperspectral (HSI) demixing \cite{HSI}, the endmembers (source signals) exhibit either mismatchs to their lab measured spectral signatures (for example, shifting along the wavenumber axis) or random distortion (some peaks are being compressed or expanded), these spectral variability are caused by many factors such as lighting conditions, geographic locations, seasons, etc.  For instance, consider a pixel (corresponding to an observation in $X$) of a mixture of soil and vegetation from a hyperspectral image, the optical reflectance spectra of soils from different geographic sites or under different lighting conditions are generally different.  Similar phenomena has been observed in Raman spectroscopy, where spectral lines of source species can have random shifts along the wavenumber axis.  In differential optical absorption spectroscopy (DOAS) \cite{DOAS1}, spectra are collected by moving a grating motor.  Changes in grating position and temperature lead to changes in dispersion of the light beam on the detector.  In nuclear magnetic resonance spectroscopy (NMR \cite{NMRBack}),  laboratory NMR spectrometers usually produce spectra of high resolution.  However, spectral variations may exist in the data as the instruments age or the changing environment.

To handle the spectral variabilities, experimentalists often adjust the reference spectra on a lookup table by knowledge of the sensing process.  In hyperspectral processing, researchers have proposed various methods to deal with the spectral variability for unmixing the observed pixels, the methods include deterministic approaches and statistical approaches (readers are referred to \cite{HSI} and references therein).  For example, authors in \cite{MoreBasis} design a selection criterion to include as many as possible variety of reflectance spectra of similar endmembers, the hyperspectral data is then unmixed with all possible different combinations of endmember signatures.  In hyperspectral target detection applications, a robust matched filter was designed in \cite{TargetDetection} by allowing the mismatch between the target source spectrum and its reference within a predefined $\epsilon$ sphere.  In \cite{L1Template}, an $\ell_1$ minimization based approach for robust template matching (data matching) in hyperspectral classification and target detection.  Other statistical approach includes Bayesian spectral mixture analysis method \cite{Bayesian} which uses the endmember signature probability distribution in the analysis for maximally capturing the spectral variability of an image with the least number of endmembers.   Although the spectral variability in hyperspectral image processing has been extensively studied, the assumptions and methods are all seeming to be limited to specific applications and could be ad hoc in many cases.  Less has been done for other spectral mixtures including Raman, NMR, and DOAS spectroscopy.  In this work, we shall assume a database containing only one reference spectrum for each source, and the objective is to improve the estimates of the mixing matrix.  This really set our work apart from these in hyperspectral unmixing.
Given the fact that the random shifts is rather small comparing to the signal length, a natural idea is to include the derivatives of the templates into the fitting basis (matrix $S$), we then have an augmented source matrix whose additional columns are the derivatives of the spectral references.  When more information is known about the shifts such as their statistical distributions, then we should include this information into the approach.  Under this assumption, we consider two scenarios both of which has random shifts in source spectral lines.  In the first case, the random shifts from source spectra are assumed to be independent and identically distributed.  In the second scenario, the random shifts are assume to have serial correlations (the shifts in mixtures from one source follow an autoregressive model), which cause correlations between mixtures.  In the first case, we can solve the problem in a sequential manner, i.e., treat one observation at a time to achieve the mixing coefficients of the source for that particular observation.  While the observations can not be treated separately if the shifts are serially correlated, they are coupled all together.

\begin{figure}
\includegraphics[height=8cm,width=15cm]{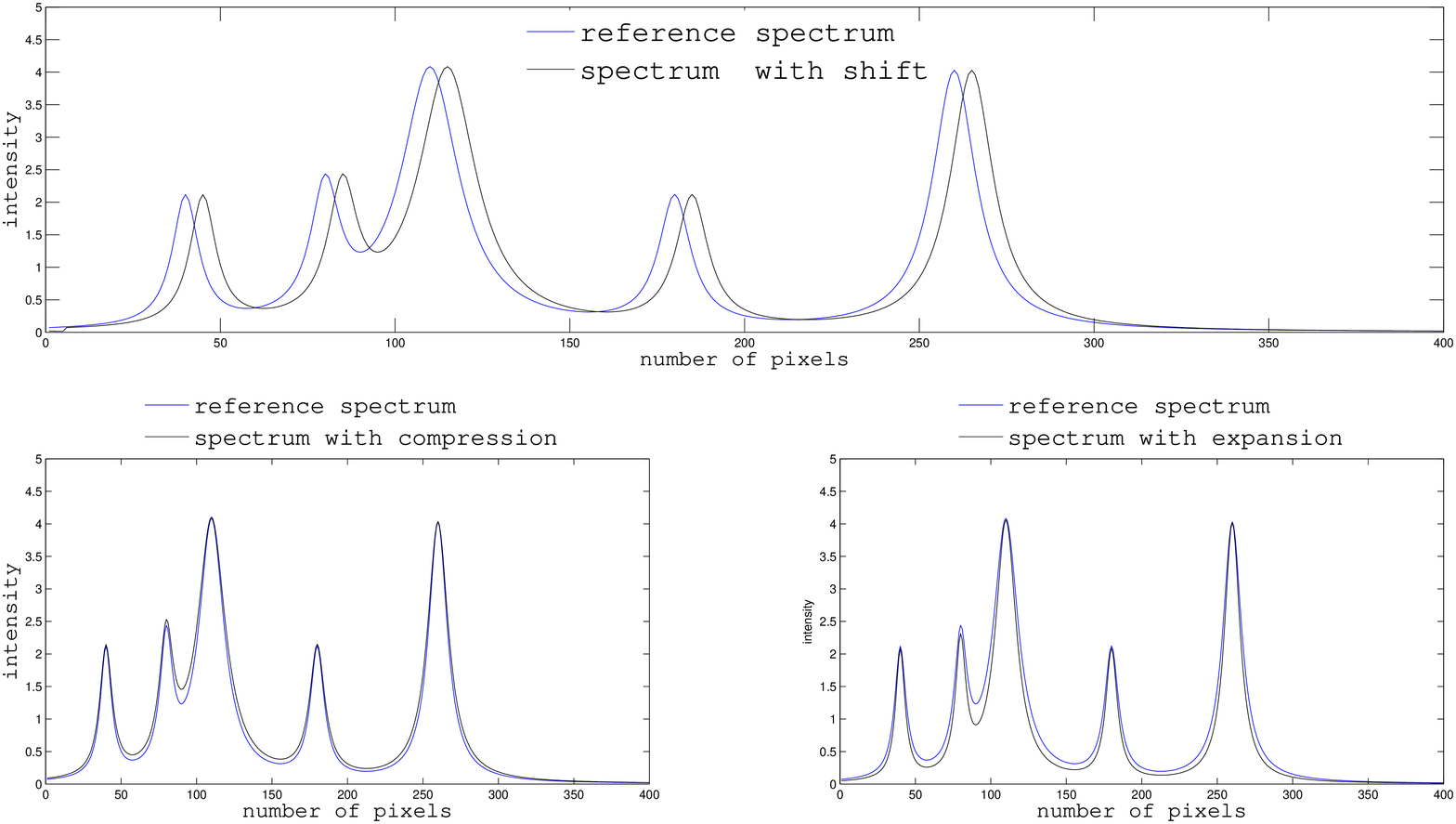}
\caption{spectral distortions.}
\label{deformation}
\end{figure}

The paper is structured as follows; In section 2, we develop the methods and their mathematical details. In section 3, numerical investigations are performed to validate the methods proposed, in addition comparisons between the methods are provided; Concluding remarks and future works are given in section 4.  {\it Throughout the paper, the following notations will be used: $\boldsymbol A$ in bold face stands for matrix, $\boldsymbol a_i$ is its ith row and $\boldsymbol a^j$ its jth column.  Greek letter $\boldsymbol \xi_{ij}$ is a random variable. }

The authors thank Professor Barbara Finlayson-Pitts for the experimental DOAS data, and Naval Research Lab for the Raman Data.  J. Xin acknowledges support of NSF grant DMS-1211179.

\section{The Methods}
 We are interested in modeling the shift and compression (or expansion) effects occurring in sensing, and build them into data matching algorithms.  Let us consider modeling shift effect by writing each row of $\boldsymbol X$ as $\boldsymbol {x}_i(\nu), i = 1,\cdots,m$, $\nu$ can be wavelength, frequency, or time depending on physical origins of the problems.  It gives the row entries when discretized.  In the model, each source signal may have different shifts in different mixtures.  Let the shift on row $\boldsymbol {s_j} = \boldsymbol {s_j}(\nu)$ of $\boldsymbol S$ be denoted by $\xi_{ij}$ which means the shift of source $\boldsymbol {s_j}$ in mixture $\boldsymbol {x_i}$.  The nonlinear mixing model with shift adjustment is:
\begin{equation}
\label{als_deformation}
{\boldsymbol x_i} = \sum^n_{j = 1}\boldsymbol{a}_{ij}{\boldsymbol s}_j(\nu + \xi_{ij}) + {\boldsymbol N}_i\;,
\end{equation}
and the related minimization problem is
\begin{equation}
\label{scaling}
(\boldsymbol {a}_{ij},\xi_{ij}) = \arg\min\|\boldsymbol x_i -\sum^n_{j = 1}\boldsymbol{a}_{ij}{\boldsymbol s}_j(\nu + \xi_{ij}) \|_2\;,
\end{equation}
 and it is non-convex. A similar problem was studied in \cite{DOAS2} based on alternating least squares estimation with Levenberg-Marquart iterative method.  In the context of image registration, the author in \cite{ImReg} propose a convex approximation of (\ref{scaling}).  Both two methods assume that each source signal has the same shift in all mixtures.  In practice, shifts of the same source may be different in different mixtures, and this is more difficult and complicated.
\subsection{Augmented Least Squares}
The idea of augmented least squares is to fit the reference spectra with their derivatives to the observations.  Given that the shifts are small comparing with the signal length, we shall use a truncated Taylor expansion of $s_j{\nu + \xi_{ij}}$ as an approximation,
\begin{equation}
{\boldsymbol s}_j(\nu + \xi_{ij}) = \boldsymbol{s}_j(\nu) + \xi_{ij} \boldsymbol{s}'_j(\nu) + \frac{1}{2} \xi^2_{ij} \boldsymbol{s}''_j(\nu) + \cdots\;,
\end{equation}
given the existence of the derivatives of $\boldsymbol {s}_i$.  In practice, it usually uses the first several terms as an approximation, in many tested cases, the first and second derivatives work well, more terms may be included as needed.  Knowing that the each mixture actually is linear combination of all the sources and their derivatives, we  include both the reference spectra and their derivatives to a form an augmented fitting basis as
$$\boldsymbol{\hat{S}} = [\boldsymbol{s_1;s'_1;s_2;s'_2;\cdots; s_n; s'_n}]$$ where only the first derivatives are included.  We will then solve the following augmented least squares problem (AgLS)
\begin{equation}
\label{eLS}
\min_{\boldsymbol{\tilde{A}}}\|\boldsymbol{X}-\boldsymbol{\tilde{A}} \boldsymbol{\tilde{S}}\,\|^2_2\;,
\end{equation}
for the linear mixture model
\begin{equation}
\label{LMM2}
\boldsymbol{X}=\boldsymbol{\tilde{A}} \,\boldsymbol{\tilde{S}} + \boldsymbol{N}
\end{equation}
The solution is $\boldsymbol{\tilde{A}} = \boldsymbol{X}\boldsymbol{\tilde{S}}^{\mathrm{T}} (\boldsymbol{\tilde{S}} \boldsymbol{\tilde{S}}^{\mathrm{T}})^{-1}$, here $\boldsymbol{\tilde{A}} $ contain weights of all the source signals and their derivatives.

Besides modeling the shift effect, augmented least squares is able to solve compression/expansion effects occurring in source signals for data matching (plots in Fig. \ref{deformation}). In fact, let the compression/expansion scale on each source signal $\boldsymbol s_j$ be denoted by $v_j$ whose value is close to 1. The nonlinear mixing model with distortion adjustment is:
\begin{equation}
\label{els_deformation}
{\boldsymbol x_i} = \sum^n_{j = 1}a_{ij}{\boldsymbol s}_j(v_j\nu) + {\boldsymbol N}_i\;, \mathrm{where}\; i = 1,\cdots,m\;,
\end{equation}
and the related minimization problem is
\begin{equation}
\label{scaling}
(\boldsymbol {A,v}) = \arg\min\|\boldsymbol x_i -\sum^n_{j = 1}a_{ij}{\boldsymbol s}_j(v_j\nu) \|_2\;.
\end{equation}
This minimization is non-convex.  Using the same idea as the augmented least squares, let $v_j = 1 + \delta_j$, then
$${\boldsymbol s}_j(v_j\nu) = {\boldsymbol s}_j(\nu + \delta_j\nu) = {\boldsymbol s}_j(\nu) + {\boldsymbol s}'_j(\nu)\delta_j\nu + \cdots. $$   Then the rest follows the equations (\ref{eLS}), (\ref{LMM2}).

  \subsection{Augmented Maximum Likelihood Estimators}

In this augmented least squares method, no prior information such as the statistical distributions, means, etc., of the random shifts are assumed.  Hence it is a rather general approach and can be applied in many situations when the knowledge of the random shifts are limited.  If certain statistics about the shifts are known a priori, we should take into account the information.

Let $\xi_{1i},\xi_{2i},\cdots,\xi_{mi}$ be the random shifts from the $ith$ source signal in all the mixtures.  We rewrite model (\ref{LMM2}) by including the source spectral and their first derivatives along with random shifts variables.
\begin{equation}
\boldsymbol{X} = \boldsymbol{A}\boldsymbol{S} + \bigl (\boldsymbol{A}\odot \boldsymbol{\Xi} \bigr )\boldsymbol{S'} + \boldsymbol{N}\;,
\end{equation}
where matrix $\boldsymbol{S'}$ contains the first derivatives of entries of $\boldsymbol{S}$, and columns of $\boldsymbol{\Xi}$ are vectors $[\xi_{1i},\xi_{2i},\cdots,\xi_{mi}]^\mathrm{T}, i = 1,\cdots,n$.  Symbol $\odot$ means the element-wise multiplication.  Below we have the model in form of matrix elements,
\begin{equation}
x_{ij} = \sum^n_{k=1}a_{ik}s_{kj} + \sum^n_{k = 1} a_{ik}\xi_{ik} s'_{kj} + n_{ij}\;.
\end{equation}
Clearly, the randomness and uncertainty of the observations include two parts: noises $\boldsymbol{N} $ and the random shift $\bigl ( \boldsymbol{\boldsymbol{A}\odot \boldsymbol{\Xi} \bigr )\boldsymbol{S'}}$.  Then the mixture matrix $\boldsymbol{X}$ has mean $\boldsymbol{AS}$ and follow the same distribution as $\bigl (\boldsymbol{A}\odot \boldsymbol{\Xi} \bigr )\boldsymbol{S'} + \boldsymbol{N}$.

Let $\Gamma_k \in \mathbb{R}^{m \times p}$ with $\Gamma_{k,ij} = a_{ik}s'_{kj}, k = 1,\cdots,n$. Notice that $\Gamma_k = \Gamma_k(A)$ depends on the unknown mixing matrix $A$, then we have
\begin{equation}
\label{LMMij}
x_{ij} = \sum^n_{k=1}a_{ik}s_{kj} + \sum^n_{k = 1} \Gamma_{k,ij} \xi_{ik} + N_{ij}
\end{equation}
Below we discuss two cases according to the correlations among the mixtures.
\subsubsection{Statistical Model Case I: Heteroscedasticity}
Assume that $\xi_{ik}$ are independent and identical Gaussian $\mathcal{N}(0,\sigma^2_k), k = 1,\cdots,n$, while $N_{ij}\thicksim \mathcal{N}(0,\tau^2)$, and $\xi_{ik} $ and $N_{ij}$ are independent. Clearly, we have the expectation of $x_{ij}$,
$\displaystyle \mathbf{E}(x_{ij}) = \sum^{n}_{k = 1} a_{ik}s_{kj}$, and the covariance between entries of $\boldsymbol X$ are
\begin{eqnarray*}
 \mathbf{Cov}(x_{ij},x_{i'j'})& = & \mathbf{Cov} (\sum^n_{k = 1}\Gamma_{k,ij}\xi_{ik},\sum^n_{k = 1}\Gamma_{k,i'j'}\xi_{ik} )\\
& = & \sum^n_{k = 1} \Gamma_{k,ij} \Gamma_{k,i'j'} \sigma^2_k\;, j \neq j'\\
& = & \sum^n_{k = 1} a^2_{ik} s'_{kj}s'_{kj'}\sigma^2_k.
\end{eqnarray*}
 And we know $\mathbf{Cov}(x_{ij},x_{i'j'}) = 0, \mathrm{if}\; i\neq i' $ since different rows of matrix $\boldsymbol X$ are uncorrelated.  However, elements within each rows are correlated, for if we rewrite the linear model (\ref{LMM2}) in terms of their rows (observations)
\begin{equation}
{\boldsymbol x_i} = \sum^n_{j = 1}a_{ij}{\boldsymbol s}_j + \sum^n_{j = 1}a_{ij}\xi_{ij}{\boldsymbol s'}_j + {\boldsymbol N}_i\;, \mathrm{where}\; i = 1,\cdots,m\;,
\end{equation}
apparently the $i$th row (observation) $\boldsymbol x_i$ is independent to
the $k$-th row $\boldsymbol x_k$ for $i\neq k$.

Let $\mathrm{vec}(\boldsymbol X)\in\mathcal{R}^{mp\times 1}$ be the vectorization of the mixture matrix $\boldsymbol {X}$.  Note that $\boldsymbol X = [\boldsymbol {x}_1^\mathrm{T},\boldsymbol{x}_2^\mathrm{T},\dots,\boldsymbol{x}_m^\mathrm{T}]^\mathrm{T}$ with $\boldsymbol{x}_i$ being a row vector of length $p$, we shall stack rows of $\boldsymbol X$ to form the $\mathrm{vec}(\boldsymbol X)$, that is
\begin{equation}
\label{vecEq}
\mathrm{vec}(\boldsymbol X) = [x_{11},x_{12},\cdots, x_{1p},\cdots,x_{m1},\cdots,x_{mp}]^\mathrm{T}.
\end{equation}

Then $\mathrm{vec}(\boldsymbol X) \thicksim \mathcal{N}_{mp}(\mathrm{vec}(AS),\boldsymbol{V})$, where $\mathcal{N}_{mp}$ denotes $mp-$variate normal distribution, and
\begin{equation*}
\boldsymbol{V} = \tau^2 \boldsymbol{I}_{mp\times mp} + \sum^{n}_{k=1} \Omega_k
\sigma^2_k\ = \left(\begin{array}{ccccc}
\boldsymbol{V_1} &  & \text{\huge O}\\
  & \ddots &   \\
\text{\Huge O} &  & \boldsymbol{V_m}
\end{array}
\right)_{mp\times mp}
\end{equation*}
 where $\boldsymbol{\Omega_k} \in \mathcal{R}^{mp\times mp} = \mathrm{diag}\biggl ((\boldsymbol \Gamma_k)_{1,:}^{\mathrm{T}}(\boldsymbol \Gamma_k)_{1,:}, \cdots,(\boldsymbol \Gamma_k)_{i,:}^{\mathrm{T}}(\boldsymbol \Gamma_k)_{i,:},\cdots,(\boldsymbol \Gamma_k)_{m,:}^{\mathrm{T}}(\boldsymbol \Gamma_k)_{m,:}  \biggr )$, here $(\boldsymbol \Gamma_k)_{i,:} \in \mathcal{R}^{1\times p} = a_{ik}{\boldsymbol s}'_{k,:} $. Hence $\boldsymbol{\Omega_k}$ is a block diagonal matrix, and $\boldsymbol{V_i}$ is $p\times p$ matrix.

The log-likelihood function is
\begin{equation}
\label{LogLike}
\mathcal{L}(A,\sigma^2,\tau^2|\boldsymbol{X,S}) = -\frac{1}{2} \ln \det({\boldsymbol{V}})
-\frac{1}{2}\Bigl (\mathrm{vec}(\boldsymbol{X-AS})^{\mathrm{T}}\boldsymbol{V}^{-1}  (\boldsymbol{X-AS}) \Bigr) -\frac{1}{2}\ln(2\pi)\;.
\end{equation}
Notice that $\boldsymbol{V} = \boldsymbol{V}(\sigma^2,A)$ is determined by $\sigma^2$ and $A$.  If $\boldsymbol{V},\sigma^2$, and $\tau^2$ are known, the maximized likelihood estimator (AgMLE) of $\boldsymbol{A}$ is the solution of $m$ generalized least squares, where $\hat{\boldsymbol{A}_{i,:}}$, the AgMLE of $\boldsymbol{A_{i,:}}$ (a row vector containing the weights of source signals in the ith mixture) for the ith mixture is
\begin{equation}
\hat{\boldsymbol{A}_{i,:}} = \arg\min \|\boldsymbol{V_i}^{-1/2}(\boldsymbol{X_{i,:} - A_i S})^{\mathrm{T}}\|^2_2\;,
\end{equation}
and it has a closed form
\begin{equation*}
\hat{\boldsymbol{A}_{i,:}} = \boldsymbol{X_{i,:} V^{-1}_i S^\mathrm{T}(SV^{-1}_i S^{\mathrm{T}})^{-1} }\;,
\end{equation*}
and $\mathbf{Var}(\hat{\boldsymbol{A}}_{i,:} ) = \boldsymbol{(SV^{-1}_i S^{\mathrm{T}})^{-1}}  $.  We propose the following iterative approach for AgMLE of $(\boldsymbol{A},\sigma^2,\tau^2)$,
\begin{itemize}
\item[1,] Start with an initial $\boldsymbol{A^{(0)}}$ obtained by the ordinary least squares,
$$
\min\|\boldsymbol{X-AS}\|^2_2 $$
whose solution is $\boldsymbol{A^{(0)} = XS^{\mathrm{T}}(SS^\mathrm{T})^{-1}}$.  Then
$\Gamma^{(0)}_{k}$ whose entry is $\Gamma^{(0)}_{k,ij}= A^{(0)}_{ik}S_{kj} $.
\item[2,] Obtain $\boldsymbol \xi = (\xi_{ik})$ by solving
$$
(\xi_{ik}) = \arg\min\|(\boldsymbol{X-A^{(0)}S})_{i,:} - \sum_k (\Gamma_k)_{i,:}\xi_{i,k} \|_2^2\;.
$$
Then $\displaystyle \sigma^{(0)}_k = \sqrt{\sum^m_{i = 1}\xi^2_{ik}}/m$,
$\displaystyle \tau^{0} = \frac{1}{\sqrt{mp}}\|(\boldsymbol{X-A^{(0)}S})_{i,:} - \sum_k (\Gamma_k)_{i,:}\xi_{i,k} \|_2 $;
 $\displaystyle \boldsymbol V^{(0)}_i = \sum^n_{k=1} \boldsymbol \Gamma_k)_{i,:}^{\mathrm{T}}(\boldsymbol \Gamma_k)_{i,:}(\sigma^{(0)}_k)^2$.
 \item [3,] Update $\boldsymbol A^{(0)}$ by
 $$
 {\boldsymbol{A}^{(1)}_{i,:}} = \arg\min \|\boldsymbol{V_i}^{-1/2}(\boldsymbol{X_{i,:} - A^{(0)}_i S})^{\mathrm{T}}\|^2_2\;. $$

 \item[4,] Set $\displaystyle \boldsymbol{A}^{(0)} \leftarrow \boldsymbol{A}^{(1)} $ and iterate step 2-3 until it converges.
\end{itemize}
The resulting variance of the estimate is given by
\begin{equation*}
\mathbf{Var}(\hat{\boldsymbol{A}}_{i,:} ) = \boldsymbol{(SV^{-1}_i S^{\mathrm{T}})^{-1}}
\end{equation*}
$95\%$ confidence interval of ${\boldsymbol{A}}_{i,:} $ is $\hat{\boldsymbol{A}}_{i,:} \pm z_{0.95} \boldsymbol{(SV^{-1}_i S^{\mathrm{T}})^{-1}}_{ii}  $, where $z_{0.95}$ is the $95\%$ percentile of a standard normal density.

\subsubsection{Statistical Model Case II: Autoregressive Model}
In this part, we are concerned with a scenario where the mixtures (observations) are correlated. The correlation might due to many factors such as the dispersion of the chemical compounds examined which is in generally a function of time, hence shift from one source existing in the observation at current time will have influence on the shift in the mixture of future experiment.  Suppose the observations are associate with time, meaning that we shall have a time series of mixtures.  Then we assume that the correlation between adjacent time the random shifts are stronger than that between shifts over longer time.  Under such condition, we propose to use autoregressive model AR(1) to characterize the random shifts exhibited in the mixtures.  We still assume an iid noise $N_{ij}\thicksim \mathcal{N}(0,\tau^2)$,$i = 1,\cdots,m, j = 1,\cdots,p$.  For any source $k = 1,\cdots,n$, we assume $\xi_{1k},\xi_{2k},\cdots,\xi_{mk}$, the shifts of source $k$ in all mixtures, are correlated according to
\begin{equation}
\xi_{ik}= \rho_k \xi_{i-1,k} + u_{ik}\;, |\rho_k|<1\;
\end{equation}
with $u_{ik}\thicksim \mathcal{N}(0,\sigma^2_k),\xi_{0,k} = 0$.  This implies that $\displaystyle \mathbf{Var}(\xi_{ik}) = \frac{\sigma^2_k}{1-\rho^2_k}, \mathbf{Corr}(\xi_{ik},\xi_{i-t,k}) = \rho^t_k$.  We also assume that $\xi_{ik}$ and $N_{ij}$ are independent, furthermore, $(\xi_{1k},\cdots,\xi_{mk}) \perp (\xi_{1k'},\cdots,\xi_{mk'})$ for $k\neq k'$.   Following the model (\ref{LMMij}), we have
\begin{eqnarray*}
\mathbf{Var}(x_{ij}) & = & \mathbf{Var}(\sum^n_{k = 1}(\Gamma_{k,ij} \xi_{ik}) + \tau^2\;,\\          & = & \sum^n_{k = 1} \Gamma^2_{k,ij}\frac{\sigma^2_k}{1-\rho^2_k} + \tau^2\;, \\
\mathbf{Cov}(x_{ij},x_{i'j'}) & = & \mathbf{Cov}(\sum^n_{k = 1}(\Gamma_{k,ij}\, \xi_{ik}, \sum^n_{k = 1} \Gamma_{k,i'j'}\, \xi_{i'k}) \\
                           & = & \sum^n_{k = 1} \mathbf{Cov}(\Gamma_{k,ij}\, \xi_{ik}, \Gamma_{k,i'j'}\, \xi_{i'k})\\
                           & = &  \sum^n_{k = 1} \Gamma_{k,ij} \Gamma_{k,i'j'}\frac{\sigma^2_k}{1-\rho^2_k}\rho^{|i-i'|}_k\;, j\neq j'\;.\\
            \end{eqnarray*}
We take the vectorization of the mixture matrix $\boldsymbol X$ as defined in (\ref{vecEq}), then $\mathrm{vec}(\boldsymbol X) \thicksim \mathcal{N}_{mp}(\mathrm{vec}(AS),\boldsymbol{V})$, and
\begin{equation*}
\boldsymbol{V} = \tau^2 \boldsymbol{I}_{mp\times mp} + \sum^{n}_{k=1} \Omega_k\frac{\sigma^2_k}{1-\rho_k}\;,
\end{equation*}
where
\begin{equation}
\boldsymbol{\Omega_k} = \left(\begin{array}{ccccc}
\boldsymbol{\Omega_{k,11}} & \boldsymbol{\Omega_{k,12}} & \ldots & \boldsymbol{\Omega_{k,1m}}\\
  \vdots & \vdots &  & \vdots   \\
\boldsymbol{\Omega_{k,m1}} & \boldsymbol{\Omega_{k,m2}}  & \ldots & \boldsymbol{\Omega_{k,mm}}
\end{array}
\right)_{mp\times mp}\;,
\end{equation}
here $\boldsymbol{\Omega_{k,ij}} = (\boldsymbol \Gamma_k)_{i,:}^{\mathrm{T}}(\boldsymbol \Gamma_k)_{j,:}\rho_k^{|i-j|}$ which is no longer a diagonal matrix.  We write
\begin{equation}
\boldsymbol{V} = \left(\begin{array}{ccccc}
\boldsymbol{V_{11}} & \boldsymbol{V_{12}} & \ldots & \boldsymbol{V_{1m}}\\
  \vdots & \vdots &  & \vdots   \\
\boldsymbol{V_{m1}} & \boldsymbol{V_{m2}}  & \ldots & \boldsymbol{V_{mm}}
\end{array}
\right)_{mp\times mp}\;
\end{equation}
where $\displaystyle \boldsymbol{V_{ij}} \in \mathbb{R}^{p\times p} = \sum^n_k (\boldsymbol \Gamma_k)_{i,:}^{\mathrm{T}}(\boldsymbol \Gamma_k)_{j,:}\rho_k^{|i-j|}\frac{\sigma^2_k}{1-\rho_k} + \tau^2 \boldsymbol{I}_{p\times p}\mathbb{I} (i = j)$, where $\mathbb{I}$ is the indicator function.

The log-likelihood function takes the same form as (\ref{LogLike}).  We propose the following iterative algorithm for estimating the parameters,
\begin{itemize}
\item[1,]
Start with an initial $\boldsymbol{A^{(0)}}$ obtained by the ordinary least squares.
\item[2,] Obtain $\boldsymbol \xi = (\xi_{ik})$ by solving
$$
\hat{\xi}_{ik} = \arg\min\|(\boldsymbol{X-A^{(0)}S})_{i,:} - \sum_k (\Gamma^{(0)}_k)_{i,:}\xi_{i,k} \|_2^2\;.
$$
Then we regress $\hat{\xi}_{ik}$ on $\hat{\xi}_{i-1,k}$ for each $k$ and we obtain
$\displaystyle \rho^{(0)}_k = \frac{\sum^m_{i=2}\hat{\xi}_{ik} \hat{\xi}_{i-1,k} }
{\sum^m_{i=2}\hat{\xi}^2_{i-1,k}}  $, and $\displaystyle \sigma^{(0)}_k = \frac{\sqrt{\sum^m_{i = 1}\xi^2_{ik}}}{m}$,
$\displaystyle \tau^{0} = \frac{1}{\sqrt{mp}}\|(\boldsymbol{X-A^{(0)}S})_{i,:} - \sum_k (\Gamma_k)_{i,:}\xi_{i,k} \|_2 $.
 \item [3,] $\displaystyle \boldsymbol V^{(0)}_{ij} = \sum^n_{k=1} \boldsymbol \Gamma_k)_{i,:}^{\mathrm{T}}(\boldsymbol \Gamma_k)_{j,:}[\rho_k^{(0)}]^{|i-j|}\frac{(\sigma^{(0)}_k)^2}{1-\rho^{(0)}_k} + (\tau^{(0)})^2 \boldsymbol{I}_{p\times p}\mathbb{I} (i = j)  $.  Then update $\boldsymbol A^{(0)}$ by
 $$
 {\boldsymbol{A}^{(1)}} = \arg\min \|\boldsymbol{V^{(0)}}^{-1/2}\mathbf{Vec}(\boldsymbol{X - A^{(0)} S})^{\mathrm{T}}\|^2_2\;. $$

 \item[4,] Set $\displaystyle \boldsymbol{A}^{(0)} \leftarrow \boldsymbol{A}^{(1)} $ and iterate step 2-3 until it converges.

\end{itemize}
\section{Experimental Investigation}
We report here the numerical examples and results of the methods. We first present the results on synthesized data, and then we show how the methods work with the realistic data from Raman, NMR, and DOAS spectroscopy.  Note that the methods require certain smoothness of the data to guarantee the continuity of the derivatives,and this condition is satisfied approximately by the test we tested.
\subsection{Synthesized examples}
\begin{figure}
\includegraphics[height=7cm,width=14cm]{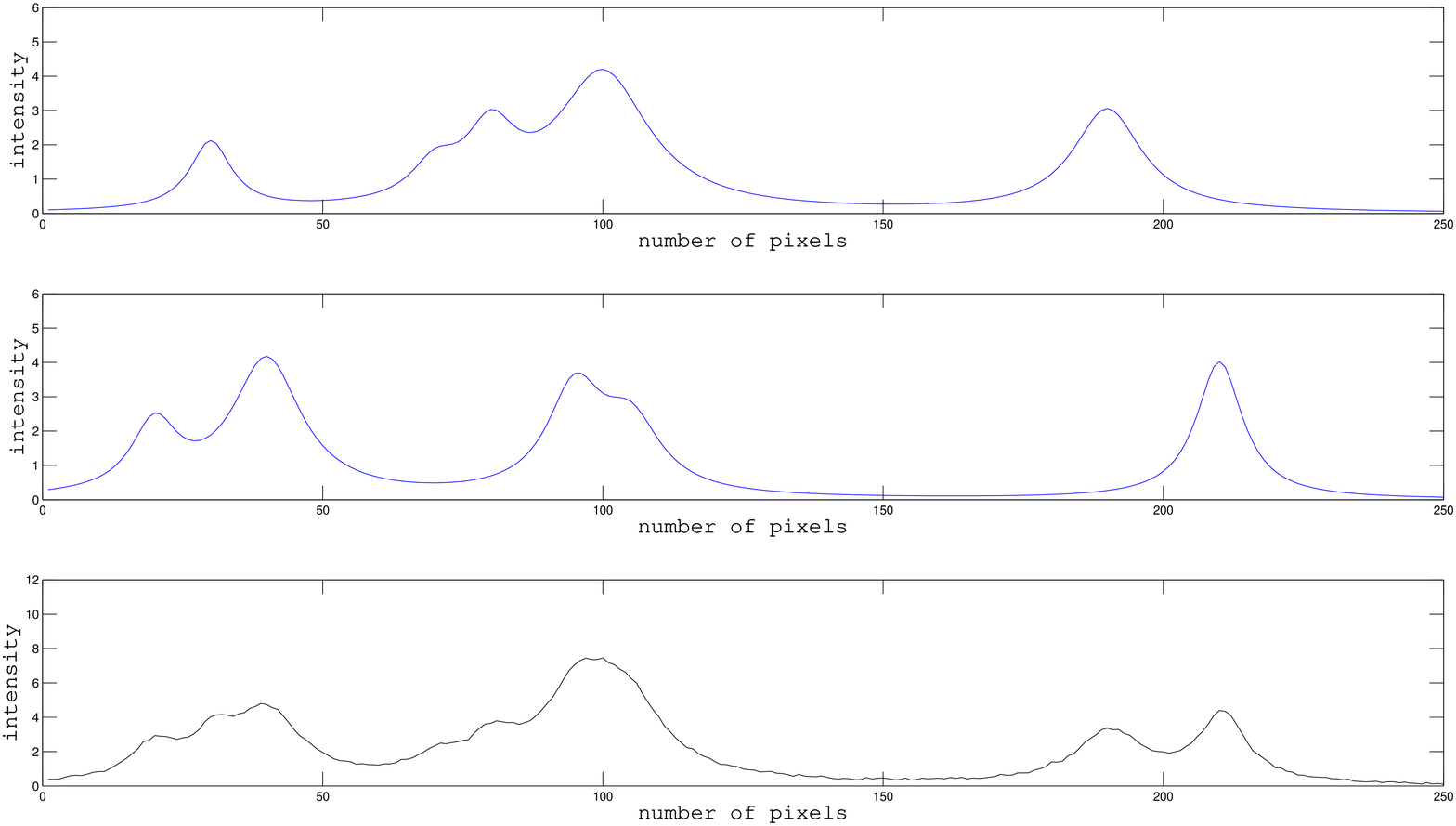}
\caption{The first two rows (blue) are the source signals, and the last row is one noisy mixture with random shift from source signals (black). }
 \label{sourceSyn}
\end{figure}

\begin{figure}
\includegraphics[height=7cm,width=14cm]{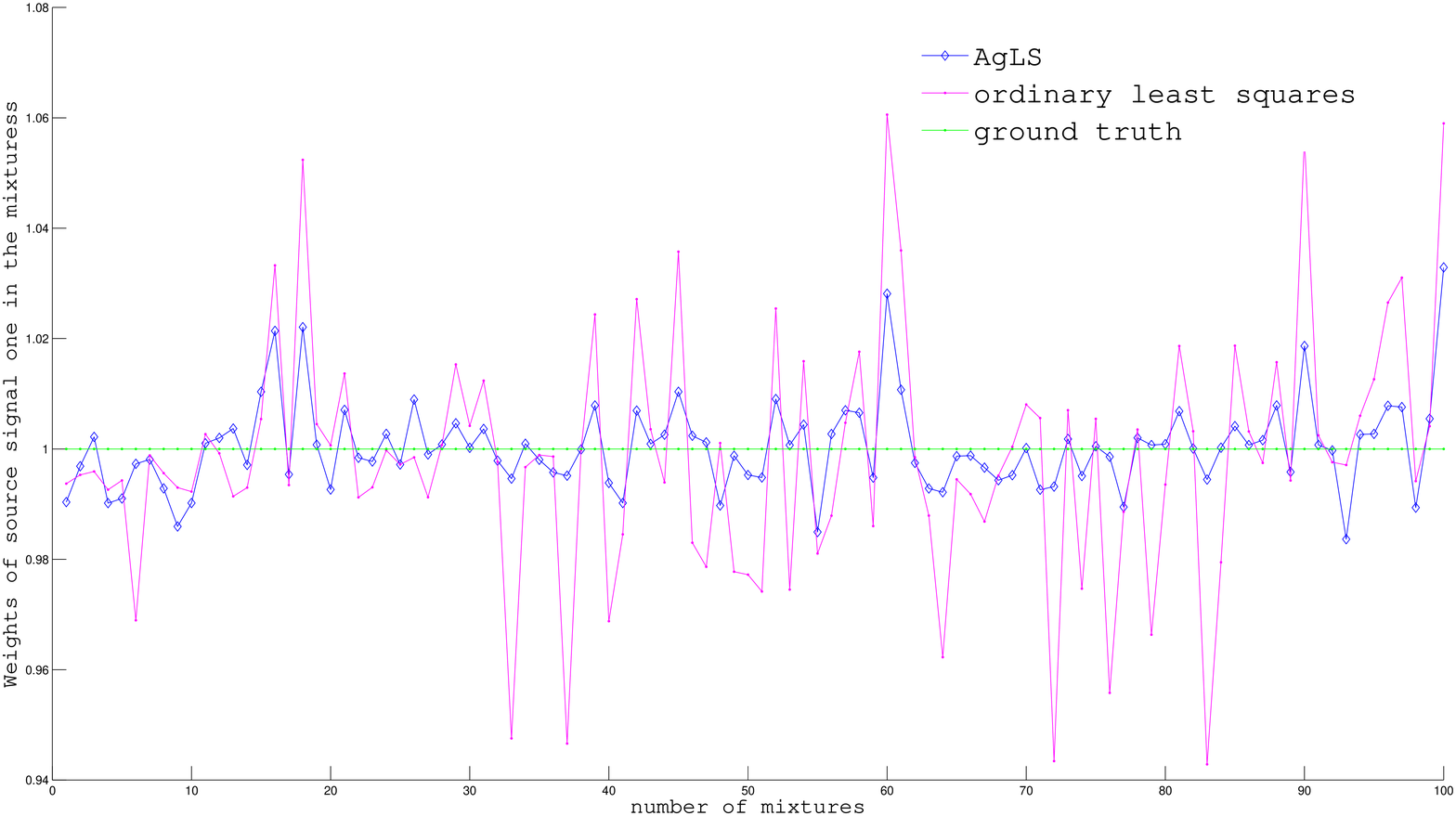}
\caption{The result of the augmented least squares (AgLS), and its comparison to the ordinary least squares. The stand deviations of the random shifts of source signals are $\sigma_1 =\sigma_2 = 1$.   The standard deviation of the noise is 0.05, the ground truth of the mixing coefficients (weights) of source 1 is one. }
\label{synResult1}
\end{figure}

\begin{figure}
\includegraphics[height=7cm,width=14cm]{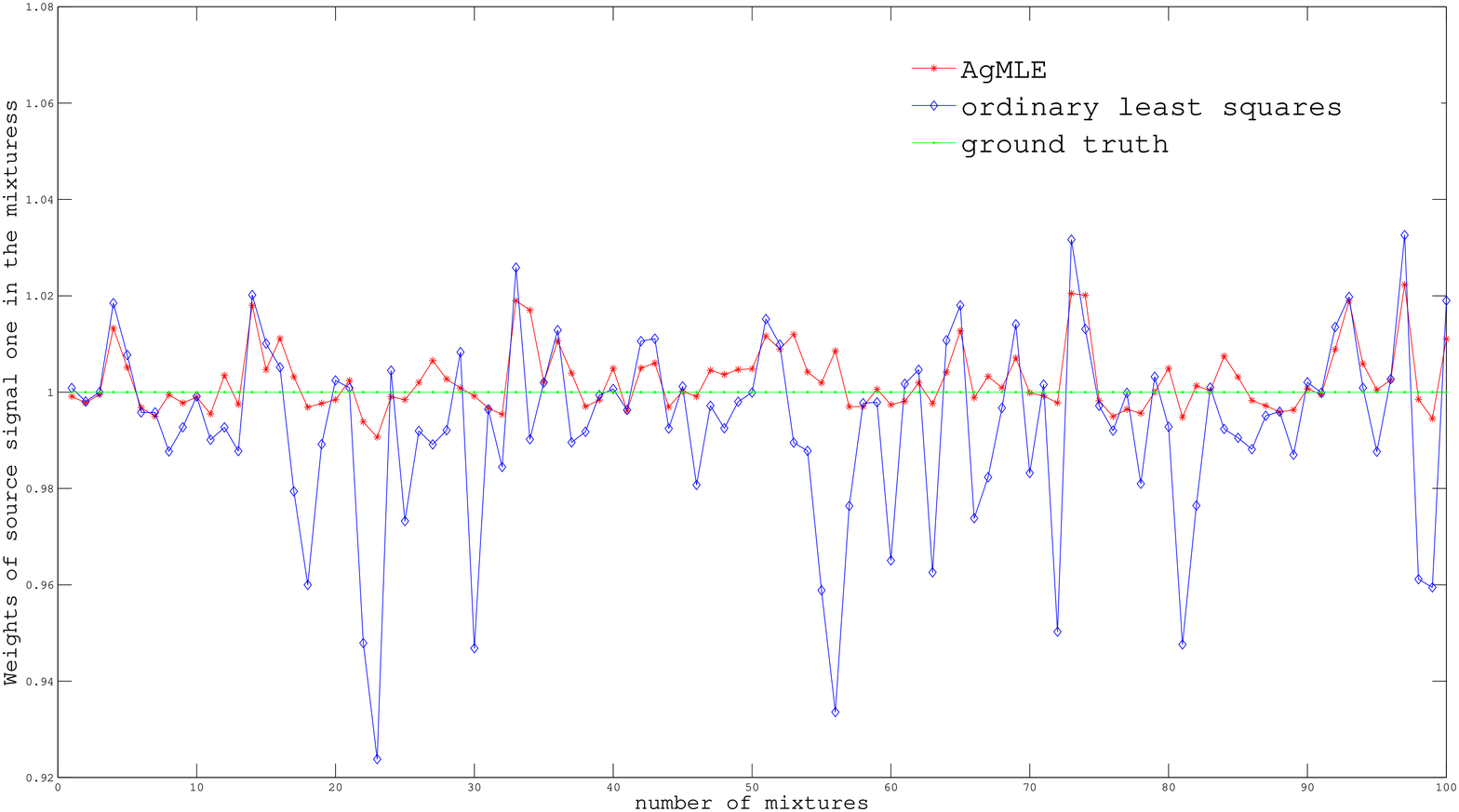}
\caption{Heteroscedasticity: the result of the augmented maximum likelihood estimator (AgMLE), and its comparison to the ordinary least squares. The stand deviations of the random shifts of source signals are $\sigma_1=1$ and $\sigma_2 = 1$, whose estimations are $\hat{\sigma}_1 = 1.0722, \hat{\sigma}_2 = 0.7617$. The standard deviation of the noise is $\tau = 0.05$ and its estimation is $\hat{\tau} = 0.1806$, the ground truth of the mixing coefficients (weights) of source 1 is one.}
\label{synResult2}
\end{figure}

\begin{figure}
\includegraphics[height=7cm,width=14cm]{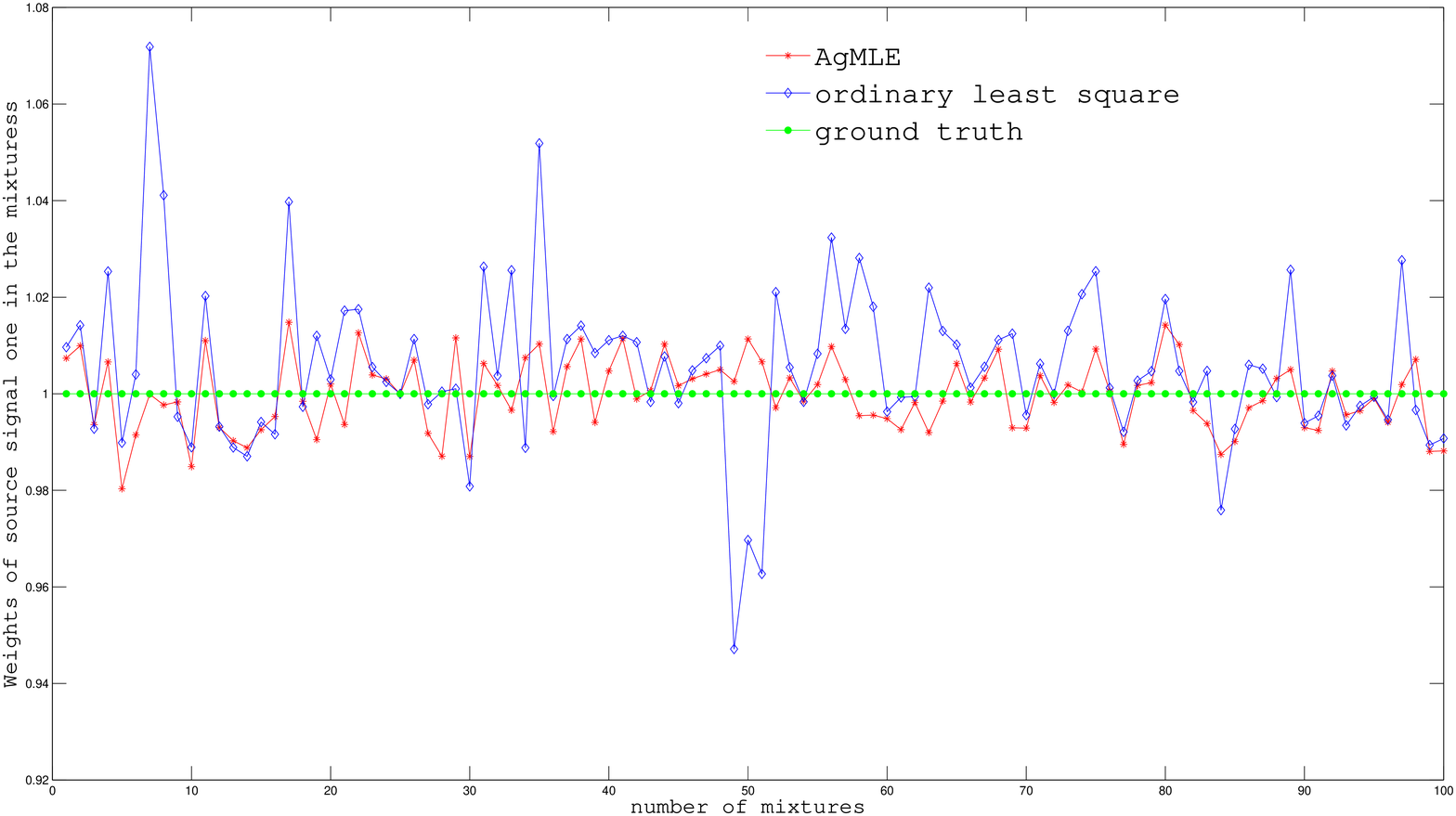}
\caption{Autoregressive model: the result of the augmented maximum likelihood estimator (AgMLE) for serial correction between mixtures, and its comparison to the ordinary least squares. The parameters used in AR(1) model are $\sigma_1 = 1, \sigma_2 = 1, \tau = 0.05, \rho_1 = 0.5, \rho_2 = 0.4$.  Their estimations are $\hat{\sigma}_1 = 0.9714, \hat{\sigma}_2 = 1.0616, \hat{\rho}_1 = 0.6168,\hat{\rho}_2 = 0.3887, \hat{\tau} = 0.2367$. The ground truth of the mixing coefficients (weights) of source 1 is one. }
\label{synResult3}
\end{figure}

\begin{figure}
\includegraphics[height=7cm,width=14cm]{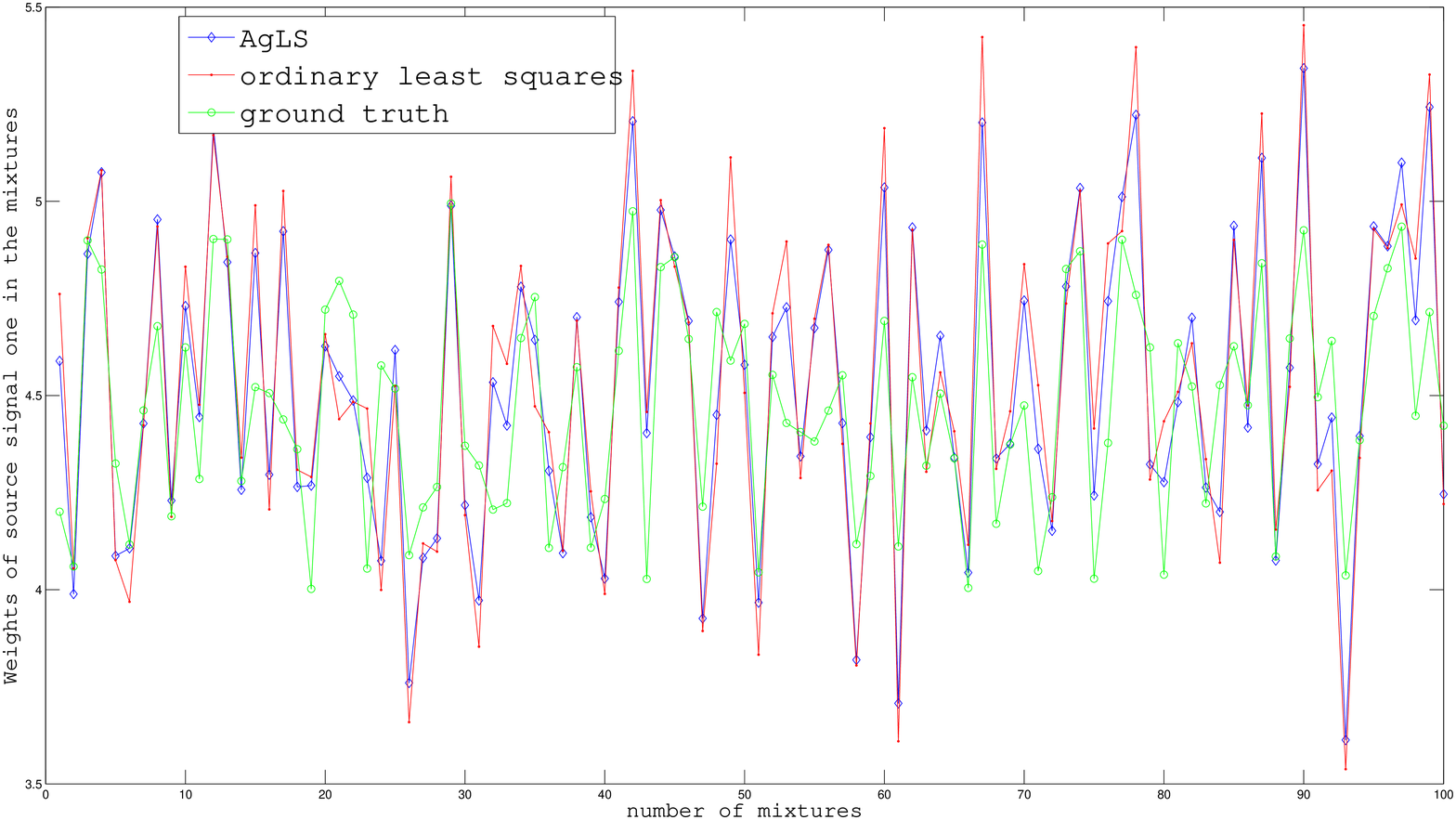}
\caption{The results of the augmented least squares (AgLS) on linear compression effect and comparisons to the ordinary least squares.  Here the compression/expansion parameters $v_j$ are in the range (0.80, 1.20), the ground truth of mixing matrix $\boldsymbol A$ are from uniform distribution of [4,5].  The absolute fitting error for the augmented least squares is 3.147, for the ordinary least squares is 3.632.}
\label{deformationELS}
\end{figure}

In this part, we show the results of simulated examples for all the three scenarios: the augmented least squares (AgLS), heteroscedasticity, and autoregressive model.  There are 2 sources, and 100 observations.  The source spectra are in Fig. \ref{sourceSyn}. The estimation of the parameters are presented below and in the Figs. \ref{synResult1}-\ref{synResult3}.  The comparisons with ordinary least squares (OLS) are also shown in the figures.  It can be seen that the result of the AgLS are better than OLS, but are less accurate than the 2 other methods when the statistics of the shifts are known.  However, the AgLS is a rather general method and computationally less complex and expensive.   We also test the effectiveness of AgLS on linear compression/expansion for signals in Fig. \ref{els_deformation}. The results are shown in Fig. \ref{deformationELS}, where the augmented least squares is able to deal with mildly linear compression or expansion.

\subsection{Applications to NMR, Raman, and DOAS spectroscopy}
\begin{figure}
\includegraphics[height=7cm,width=14cm]{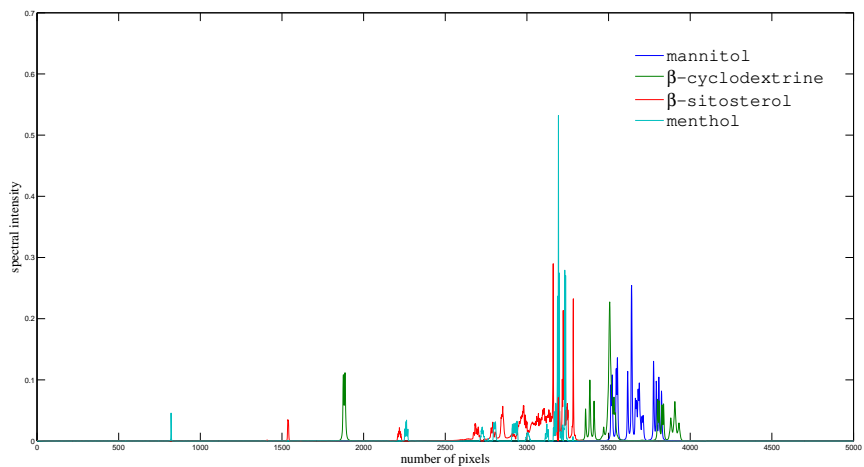}
\caption{The reference spectra of mannitol, $\beta-$cyclodextrine, $\beta-$sitosterol, and menthol. }
\label{nmrSource}
\end{figure}

\begin{figure}
\includegraphics[height=7cm,width=14cm]{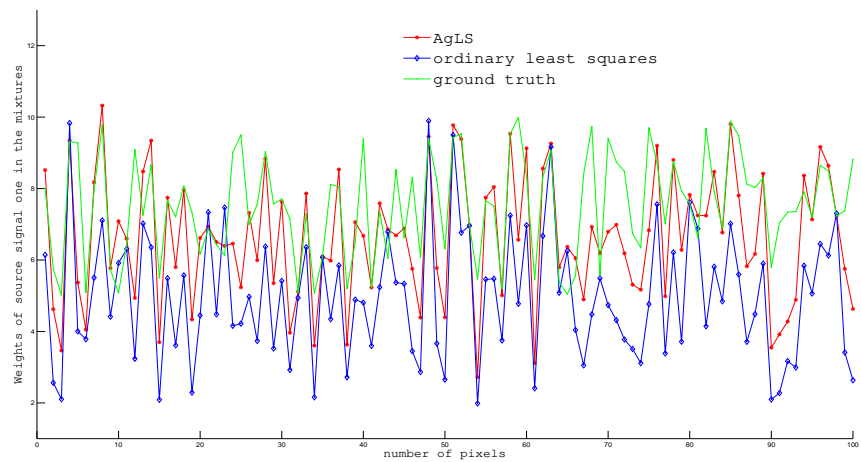}
\caption{Augmented least squares: the estimation of the weights of mannitol in the 100 mixtures. }
\label{NMR1}
\end{figure}

\begin{figure}
\includegraphics[height=7cm,width=14cm]{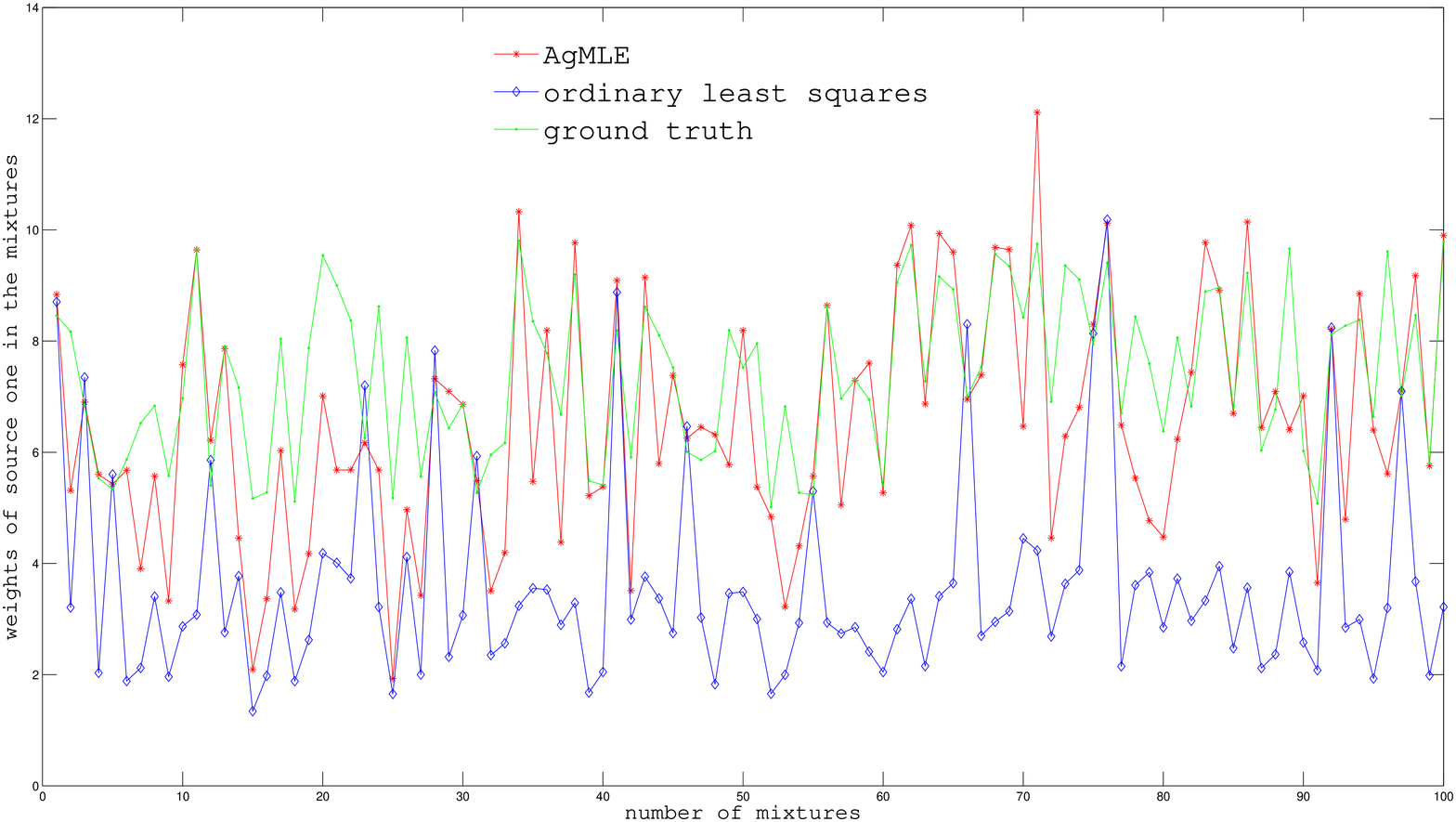}
\caption{Augmented maximum likelihood estimator: the estimation of the weights of mannitol in the 100 mixtures.  The parameters used in the example: $\sigma_i  = 2, i = 1,\cdots, 4, \tau = 0.05$.  The computed results of the parameters are $\hat{\sigma}_1 =  0.9353, \hat{\sigma}_2 = 1.0264,   \hat{\sigma}_3 = 0.9229,   \hat{\sigma}_4 =  1.0766 ,\hat{\tau} = 0.1940$. }
\label{NMR2}
\end{figure}

{\bf 1, Nuclear Magnetic Resonance (NMR) spectroscopy} is a powerful and popular tool for chemists and biochemists to investigate and determine the structures and properties of molecules.  The NMR spectrum of a chemical compound is produced by the Fourier transformation of a time-domain signal which is a sum of
sine functions with exponentially decaying envelopes.  The real part of the spectrum can be presented as the sum of symmetrical, positively valued, Lorentzian-shaped peaks.    In the example of NMR spectroscopy, we use the true spectra of four chemical compounds, mannitol, $\beta-$cyclodextrine, $\beta-$sitosterol, and menthol (see Fig. \ref{nmrSource} for their spectra). The coefficients of the mixing matrix $A$ are generated from a uniform distribution in the range of [5,10] for each mixture.  Then we obtain the mixtures by simulating the mixing process in (\ref{LMM}) with randomly shifted source spectra and by adding Gaussian noise. The results of the estimations are listed in Figs. \ref{NMR1},\ref{NMR2}.  It once again shows that the methods proposed are better than the ordinary least squares.

\begin{figure}
\includegraphics[height=7cm,width=14cm]{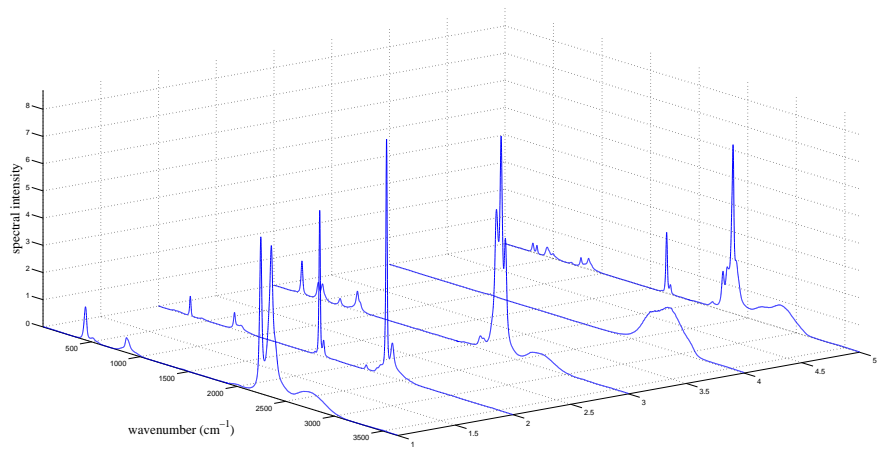}
\caption{Raman spectra of some liquids used to make explosives: from left to right, Methanol, Acetonitrile, Ethanol, Water, and their mixture. }
\label{Raman}
\end{figure}

\begin{figure}
\includegraphics[height=7cm,width=14cm]{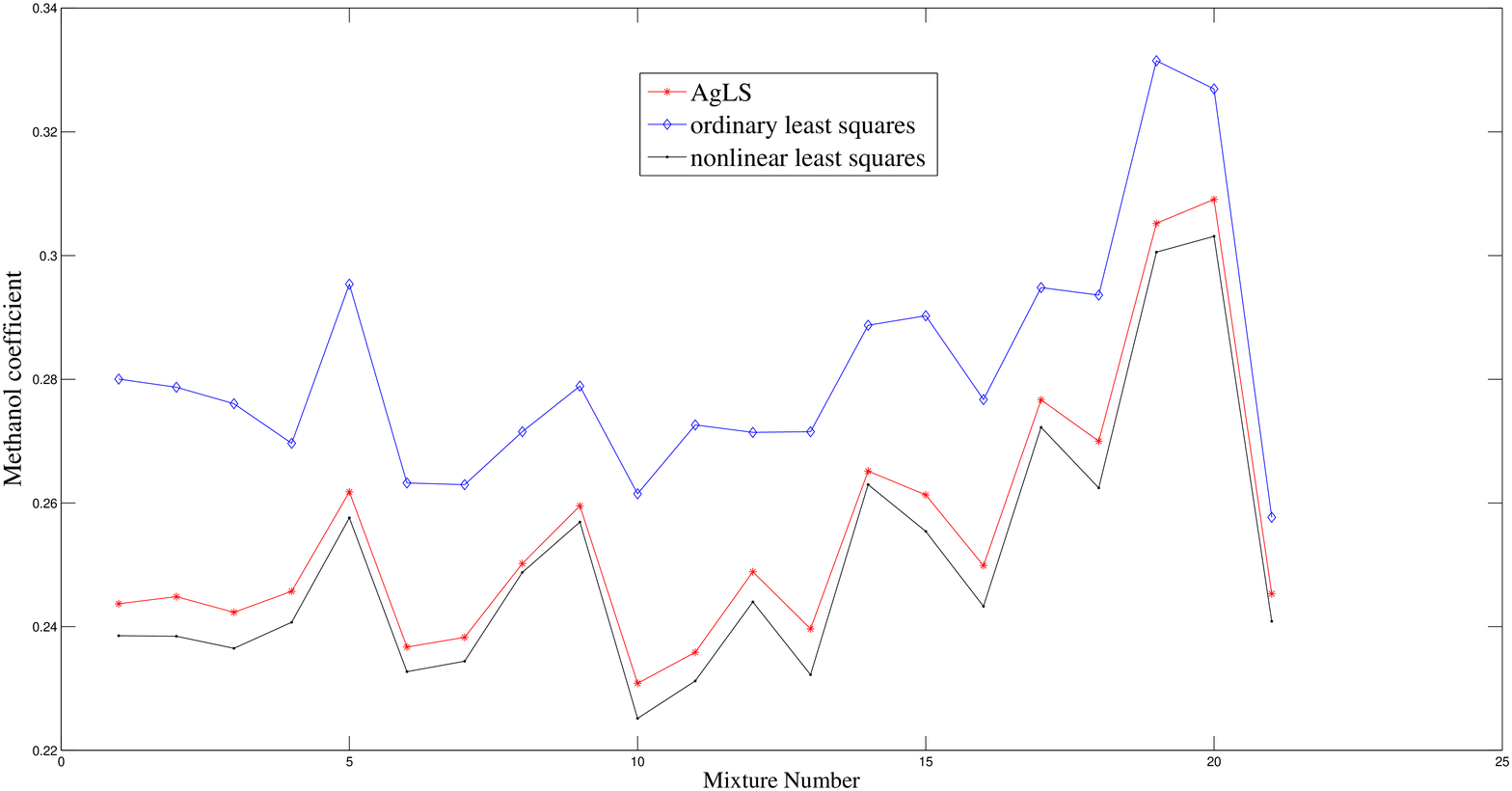}
\caption{The estimation of mixing coefficients of Methanols computed by augmented least squares.  As a comparison, we also present the results by nonlinear least squares. }
\label{RamanResult1}
\end{figure}

\begin{figure}
\includegraphics[height=7cm,width=14cm]{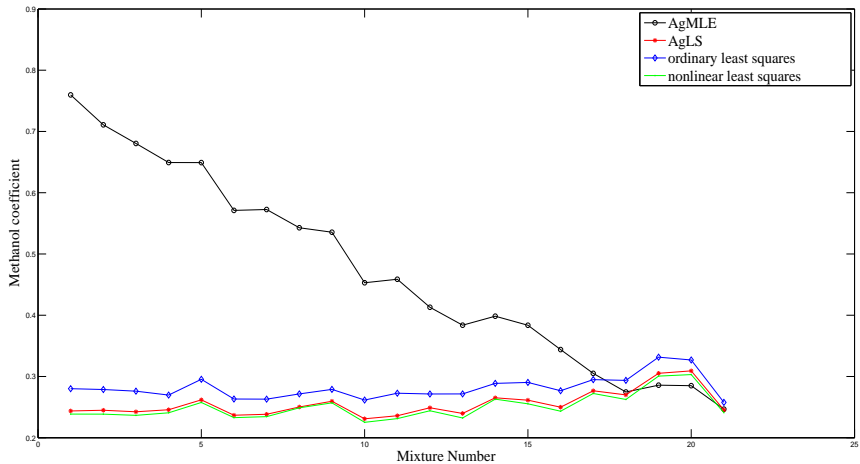}
\caption{The estimation of mixing coefficients of Methanols computed by augmented maximum likelihood (AgMLE) estimator and comparison to other methods.}
\label{RamanResult2}
\end{figure}

{\bf 2, Raman spectroscopy}:  In a second realistic example, we apply our method to the Raman spectra.  As shown in Fig. \ref{Raman}, a Raman spectrum gives a collection of peaks that correspond to the characteristic vibrational frequencies of the materials being examine, thus providing unique information on the molecular structure and chemical composition of the matters \cite{RamanBack}.  We test our method on Raman data provided by NSF/ATD program \cite{atd}.  The sample consists of several liquid substances, some of them are commonly used for making explosives.  The dataset include 21 mixed Raman spectra at different incident laser wavelengths.   These are realistic data from Naval Research Lab and we do not have the ground truth of the mixing coefficients to compare to, moreover we do not have knowledge of the statistics of the random shifts.  In this situation, it is clearly the augmented least squares is the most suitable method for estimation of the mixture coefficients.  As a comparison we also shown the comparisons to the case 2 where IID Gaussian distributed shifts are assumed. We also compare our methods to the nonlinear least solver from \cite{DOAS2}  The results of augmented least squares are close to those of nonlinear least squares (see Fig. \ref{RamanResult1}), although the augmented least squares is less complicated and easier to implement.  On the other hand, the maximum likelihood results are much more deviated from other methods for the mixtures number 1-15 (Fig. \ref{RamanResult2}).  This observation shows that shifts are not Gaussian or the existence of other nonlinear distortion such as compression or expansion in the signals.

\begin{figure}
\includegraphics[height=7cm,width=14cm]{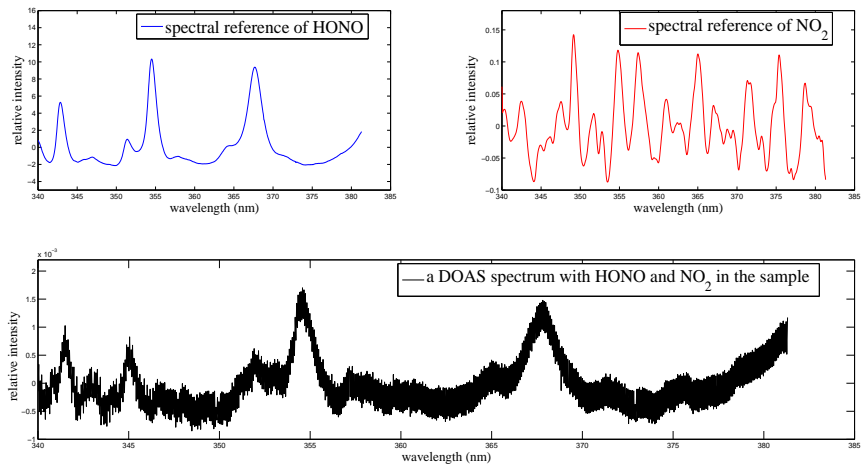}
\caption{The spectral references of trace gases HONO and $\mathrm{NO_2}$, and their noisy mixture. }
\label{DOAS}
\end{figure}

{\bf 3, Differential optical absorption spectroscopy} (DOAS) is based on the light absorption property of matter to identify broadband and narrow band spectral structures, and analyze atmospheric trace gases concentrations \cite{DOAS1}.  It can also help to understand the influence of atmospheric chemistry on climate and air quality.  Fig. \ref{DOAS} shows spectra reference of trace gases HONO and $\mathrm{NO_2}$ and a DOAS spectrum of mixture sample containing these two gases.  We are here to test the methods on realistic DOAS data.  We are glad to find that the augmented least squares is able to deliver better results compare to the ordinary least squares, on the other hand, the maximum likelihood estimator (Heteroscedasticity) produces same results as the ordinary least squares, in fact, the matrix $\boldsymbol Q$ is actually diagonally dominant, and the shift effect is insignificant.  However, the better results of augmented least squares implies that there are other spectral distortions such as linear compression or expansion.  This once again shows that the augmented least squares are able to handle not only the shift but also other distortions in the spectra.

\begin{figure}
\includegraphics[height=7cm,width=14cm]{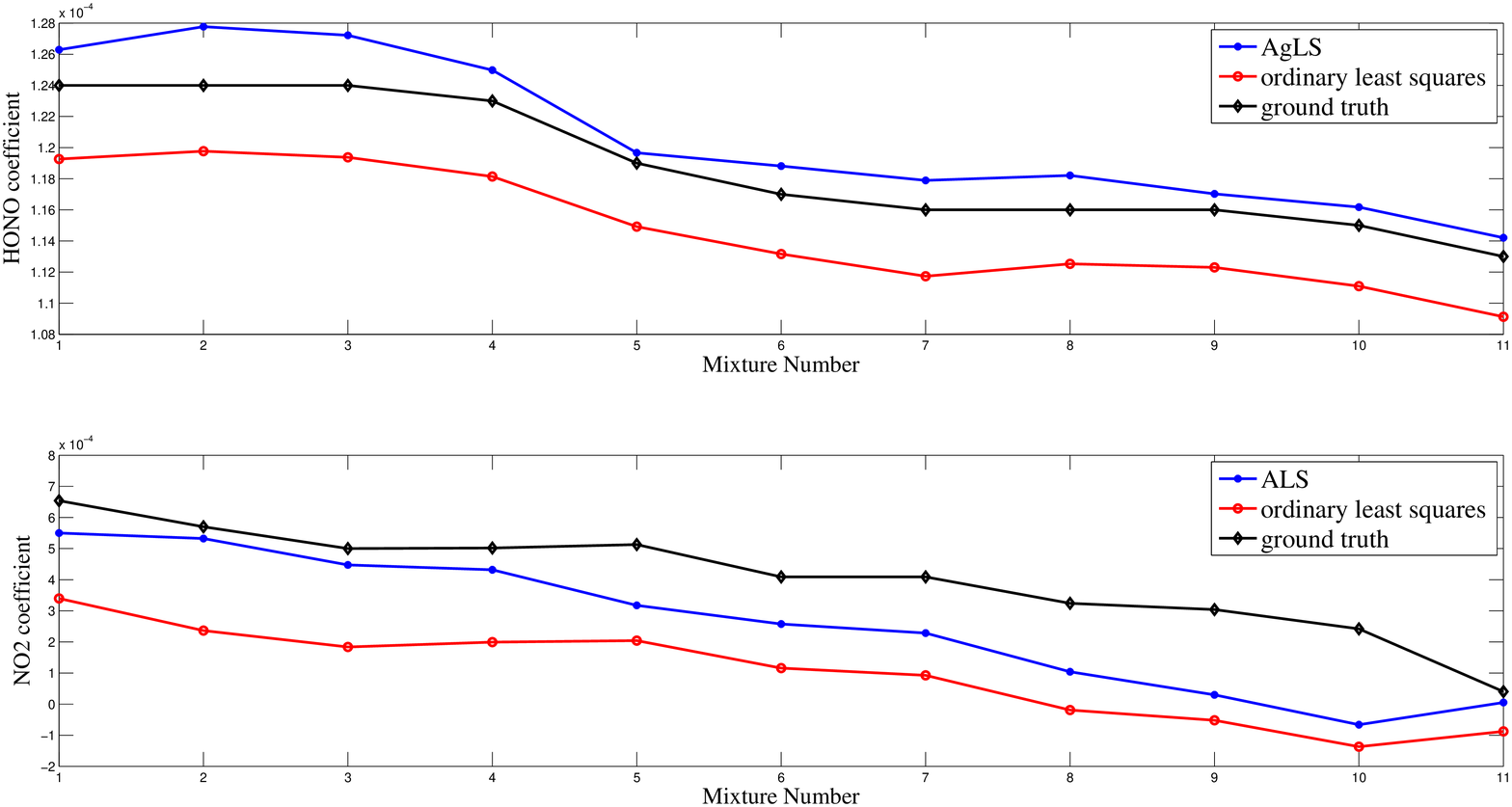}
\caption{The computed coefficients of HONO and $\mathrm{NO_2}$ by augmented least squares (AgLS) and comparisons with other methods.}
\label{DOAS_ELS}
\end{figure}

\begin{figure}
\includegraphics[height=7cm,width=14cm]{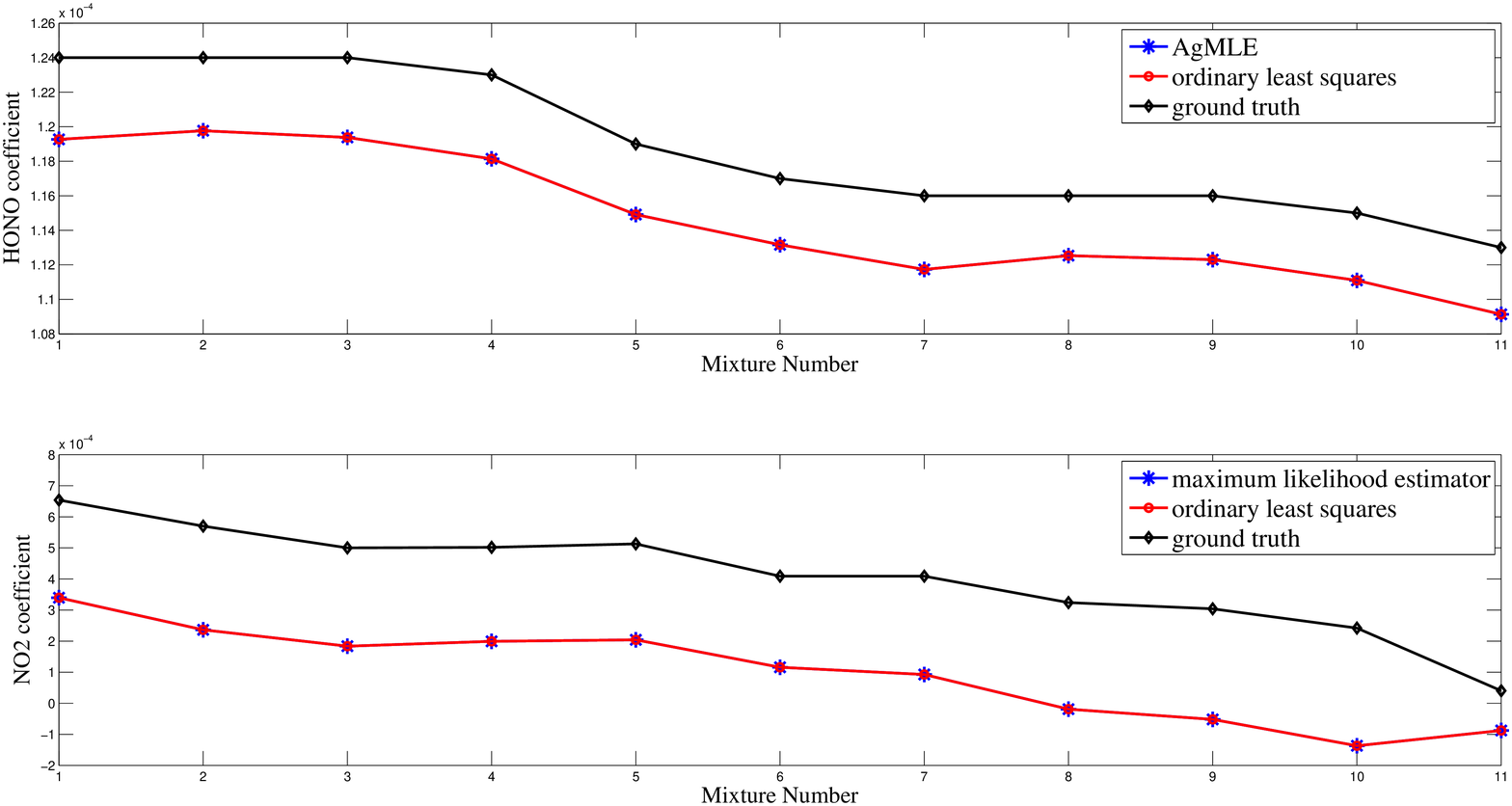}
\caption{The coefficients of HONO and $\mathrm{NO_2}$ by augmented maximum likelihood estimator (AgMLE) and comparisons to other methods.}
\label{DOAS_MLE}
\end{figure}

  \section{Conclusion and Future Works}
\begin{figure}
\includegraphics[height=3.5cm,width=7.5cm]{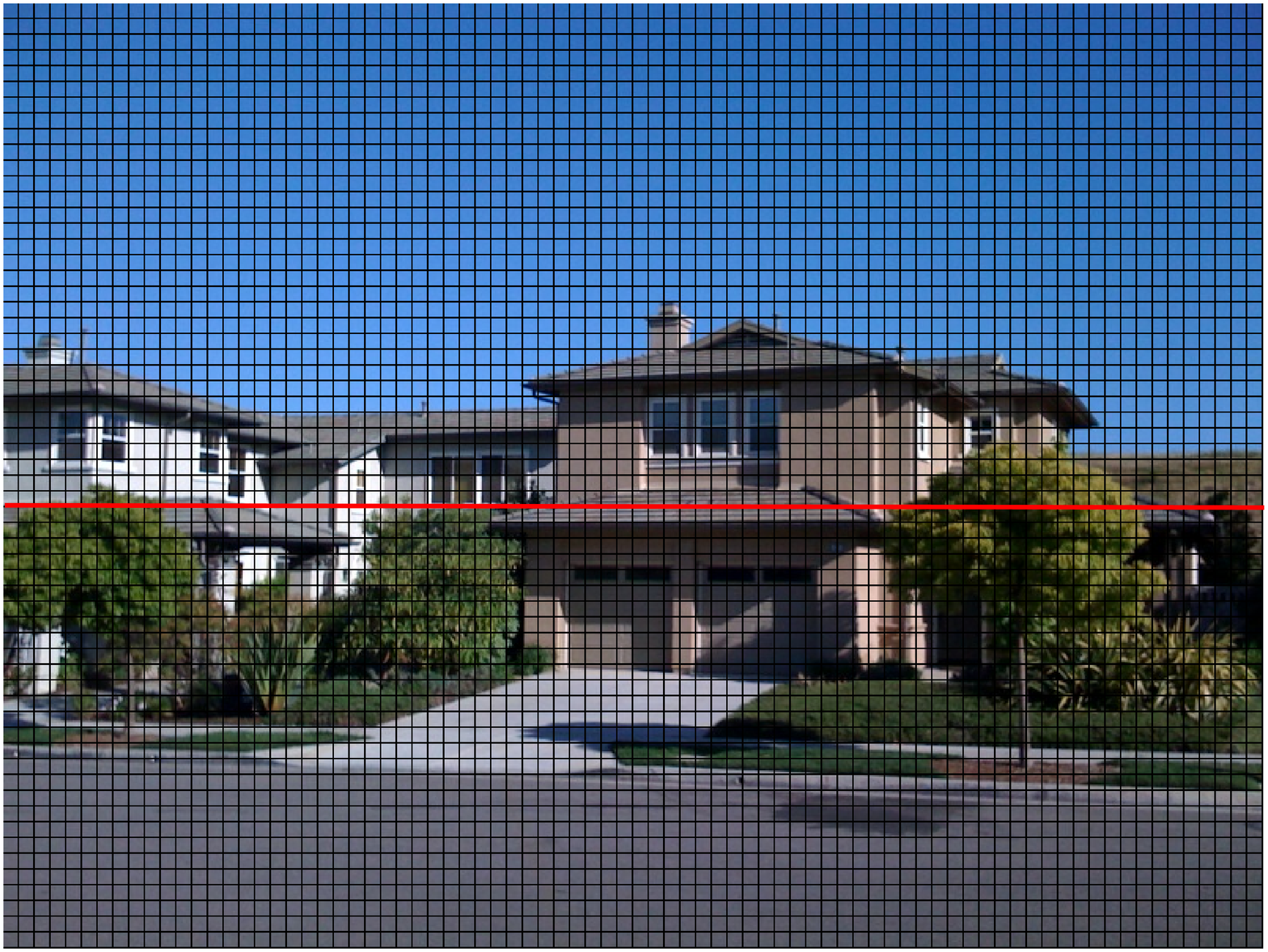}
\includegraphics[height=3.5cm,width=7.5cm]{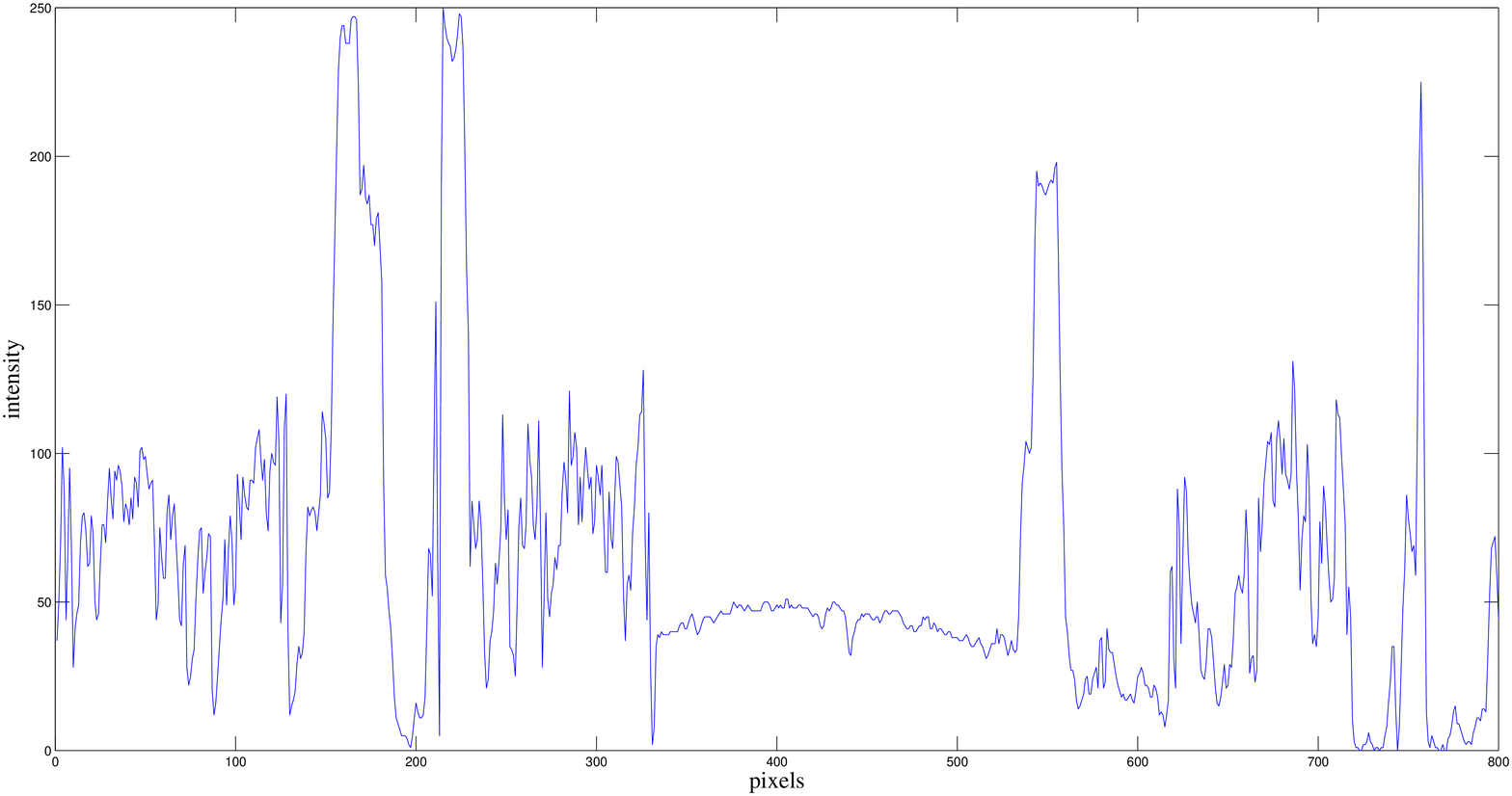}
\caption{signal with jumps and oscillations (right) is a slice (red line) of an image containing man-made buildings and trees.}
 \label{discont}
\end{figure}

In this paper, we are concerned with data matching where the source signals have random shifts and/or other nonlinear distortions.  Novel methods are proposed for various scenarios according to the knowledge of the data, and they prove to work well with real-world data from NMR and Raman as well DOAS spectroscopy.  The modeling of the nonlinear distortions is based on the truncated Taylor expansion of the signals.   The computational approach is an augmented least squares which fits the reference spectra of the signals and their derivatives to the mixtures.  If the random distortions are Gaussian, an augmented maximum likelihood estimator (AgMLE) is developed.  In future work, we shall study non-Gaussian priors (such as Laplacian or hyper Laplacian) in maximum likelihood estimator to model the shifts and other distortions.  For signals with jumps and/or high oscillations, the methods based on derivatives would fail to work.  These signals can be found in intensity image of a building with trees and vegetation for image registrations (see Fig. \ref{discont}).  One idea is to model shifting by shift operator that is characterized by sparse matrix of one and zeros.  The unknowns to be estimated are the mixing matrix and the shift matrix, and they are multiplied together, hence the problem is non-convex.  A future work is to study convex approximation of the problems.

\begin{thebibliography}{9999999}
\bibitem{atd} ATD data, http://grassmann.math.colostate.edu/ATD/home.html.

\bibitem{Bayesian} C. Song, {\em Spectral mixture analysis for subpixel vegetation fractions in the urban environment: How to incorporate endmember variability?}, Remote Sensing Environ., vol. 95, no. 2, pp. 248--263, 2005.

\bibitem{Wolke88}R. Wolke and H. Schwetlick, {\em Iteratively Reweighted Least Squares: Algorithms, Convergence Analysis, and Numerical Comparisons}, SIAM J. Sci. and Stat. Comput., 9(5), 907 -969, 1988.
\bibitem{HSI} A. Zare and K. C. Ho, {\em Endmember variability in hyperspectral analysis}, IEEE Signal Processing Mag., vol. 31, no. 1, pp. 95--104, 2014.

\bibitem{HSI1} C-I Chang, ed., {\em Hyperspectral Data Exploitation: Theory and Applications}, Wiley-InterScience, 2007.

\bibitem{ImReg} E. Esser, {\em A Convex Model of Image Registration}, UCLA CAM Report, 10-04, 2010.


\bibitem{LM} D. Marquardt. {\em An Algorithm for Least-Squares Estimation of Nonlinear Parameters}, SIAM Journal on Applied Mathematics, 11 (2), 431–441, 1963.

\bibitem{TargetDetection} D. Manolakis, R. Lockwood, T. Cooley, and J. Jacobson, {\em Robust matched filters for target detection in hyperspectral imaging data}, in Acoustics, Speech and Signal Processing, 2007. ICASSP 2007. IEEE International Conference, Vol. 1, april 2007, pp. I-529-I-532.
\bibitem{DOAS1} U. Platt and J. Stutz, {\em Differential Optical Absorption Spectroscopy: Principles and Applications}, Springer, 2008.
 \bibitem{DOAS2} J. Stutz and U. Platt, {\em Numerical analysis and estimation of the statistical error of differential optical absorption
spectroscopy measurements with least-squares
methods. Appl. Optics.,} {\bf 35}, pp. 6041--6053, 1996.

\bibitem{MoreBasis} D.A. Roberts, M. Gardner, R. Church, S. Ustin, G. Scheer, R.O. Green, {\em Mapping Chaparral in the Santa Monica Mountains using multiple endmember spectral mixture models}, Remote Sensing of Environment, 65, pp. 267--279, 1998.

\bibitem{NMRBack} R. M. Silverstein, F. X. Webster, and D. J. Kiemle,
{\em Spectrometric Identification of Organic Compounds
,} John Wiley\& Sons, 2005.
\bibitem{RamanBack} D.A. Long, {\em The Raman Effect: A Unified Treatment of the Theory of Raman Scattering by Molecules}, John Wiley \& Sons, 2002.
\bibitem{L1Template} Z. Guo and S. Osher, {\em Template Matching via $\ell_1$ Minimization and Its Application to Hyperspectral Data,} Inverse Problem and Imaging, Vol 5, No. 1, 2011, pp. 19--35.

\end{thebibliography}
\end{document}